\documentstyle[12pt,twoside]{article}
\newtheorem{theorem}{\bf Theorem}[section]

\newtheorem{lemma}[theorem]{\bf Lemma}
\newtheorem{corollary}[theorem]{\bf Corollary}

\input{amssym}

\date{}
\begin{document}

\title{{\Large\bf Jordan derivations on the $\theta-$Lau products of Banach algebras}}

\author{{\normalsize\sc M. Ghasemi and M. J. Mehdipour}}

\maketitle

{\footnotesize  {\bf Abstract.}  In this paper, we study Jordan derivation-like maps on the $\theta-$Lau products of algebras. We characterize them and prove that under certain condition any Jordan derivation-like maps on the $\theta-$Lau products is a derivation-like map. Moreover, we investigate the concept of centralizing for Jordan derivation-like maps on the $\theta-$Lau products of algebras.\\


}

\section{\normalsize\bf Introduction}

Let $A$ be a Banach algebra. Let us recall that a linear mapping $\frak{D}: A\rightarrow A$ is called a \emph{derivation} if
$$\frak{D}(ax)=\frak{D}(a)x+a\frak{D}(x)$$
for all $a, x\in A$. Also, $\frak{D}$ is called a \emph{Jordan derivation} if for every $a\in A$
$$
\frak{D}(a^2)=\frak{D}(a)a+a\frak{D}(a).
$$
The set of all derivations and Jordan derivations on $A$ are denoted by $\mbox{Der}(A)$ and $\mbox{Der}_J(A)$, respectively.

Let $B$ be a Banach algebra and $\theta$ be a nonzero multiplicative linear functional on $B$. Following [16], the $\theta-$Lau product $A$ and $B$ is denoted by $A\times_\theta B$ and it is the direct product $A\times B$ together with the component
wise addition and the multiplication
$$
(a, b)\cdot_\theta(x, y)=(ax+\theta(y)a +\theta(b)x,by).
$$
We note that in the case where $B={\Bbb C}$ and $\theta$ is the identity map on ${\Bbb C}$, the unitization $A$ will be obtained. We also note that if we permit $\theta=0$, the $\theta-$Lau product $A\times_\theta B$ is the usual direct product. Hence we disregard the possibility that $\theta=0$.

The $\theta-$Lau products $A\times_{\theta}B$ were first introduced by Lau [12], for certain Banach algebras.
Sanjani Monfared [16] extended this product to arbitrary Banach algebras $A$ and $B$. The $\theta-$Lau products are significance and utility. Because, the $\theta-$Lau product is a strongly splitting Banach algebra extension of $B$ by $A$; for the study of extensions of Banach algebras see [3, 7]. Also, many properties are not shared by arbitrary strongly
splitting extensions, while the $\theta-$Lau products exhibit them; see  [16]. Furthermore, the
$\theta-$Lau products are a source of examples or counterexamples; see for instance [17]. These
reasons caused that several authors studied various aspects of the products [6, 8, 9, 11, 17, 19]. In this paper, we continue these investigations and study Jordan derivation-like maps of them.

It is clear that every derivation is a Jordan derivation. But, the converse is, in general, not true. Here a question arises: when dose the converse hold? In 1957, Herstein [10] proved that every Jordan
derivation on a 2-torsion free prime ring is a derivation; see also [1, 5, 14, 15]. Bresar [4] gave a generalization of
Herstein's result for semiprime rings. Many attempts were made to study this question for Jordan derivations on Banach algebras [2, 4, 18]. For example, Sinclair [18] proved that every continuous Joradan derivation on a semisimple Banach algebra is a derivation. Brasar [4] showed that any Jordan derivation on a semisimple Banach algebra is continuous. So any Jordan derivation on a semisimple Banach algebra is a derivation. It is natural to ask whether results concerning Jordan derivations on Banach algebras hold for the $\theta-$Lau products $A\times_\theta B$? The other question comes to mind immediately: what happens to $\theta$ in these investigations? To answer these questions, we consider linear mappings $d: A\times B\rightarrow A\times B$ satisfying
$$
d((a, b)\cdot_{\theta}(a, b))=d(a, b)\cdot_\phi(a,b)+(a, b)\cdot_\gamma d(a, b)\,\,\,\,(a\in A\,\,\mbox{and}\,\,b\in B),
$$
where $\theta, \phi$ and $\gamma$ are nonzero multiplicative linear functional on $B$ .
We denote the set of all theses mappings by $\mbox{Der}_J(A\times_\theta^{\phi, \gamma} B)$. In this paper, we investigate the questions concerning Jordan derivations for elements of $\mbox{Der}_J(A\times_\theta^{\phi, \gamma} B)$.

In this paper, we characterize elements of $\mbox{Der}_J(A\times_\theta^{\phi, \gamma} B)$ in the case where $A$ has a right identity. We also give a necessary and sufficient condition under which every element of $\mbox{Der}_J(A\times_\theta^{\phi, \gamma} B)$ is a derivation. For unitary algebra $A$ and semisimple Banach algebra $B$, we prove that if $\theta\neq\phi$, then $\mbox{Der}_J(A\times_\theta^{\phi, \phi} B)=\mbox{Der}(A\times_\theta^{\phi, \phi} B)$. Furthermore, we investigate $(\eta_1, \eta_2)-$centralizing element of $\mbox{Der}_J(A\times_{\theta}^ {\phi, \gamma}B)$.

\section{Main Results}
In the sequel, let $A$ be a Banach algebra with a right identity $u$ and right annihilator $\mbox{ran}(A)$, the set of all $z\in A$ with $az=0$ for all $a\in A$. Let also $\theta, \phi$ and $\gamma$ be nonzero multiplicative linear functionals on any Banach algebra $B$. The next Lemma is needed to prove our results.

\begin{lemma} Let $d: A\times B\rightarrow
A\times B$ be a mapping with
$$
d((a, b)\cdot_{\theta}(a, b))=d(a, b)\cdot_\phi(a, b)+(a, b)\cdot_\gamma d(a, b)
$$
for all $a\in A$ and $b\in B$. Then $d$ maps $A$ into itself and $d(u, 0)\in\emph{\mbox{ran}}(A)$.
\end{lemma}
{\it Proof.} By hypothesis, we have
\begin{eqnarray*}
d((a+u, 0)\cdot_{\theta}(a+u, 0))&=&d(a+u, 0)\cdot_\phi(a+u, 0)\nonumber\\
&+&(a+u, 0)\cdot_\gamma d(a+u, 0)\nonumber
\end{eqnarray*}
for all $a\in A$. So
\begin{eqnarray}\label{1ud}
d(a, 0)&+&d(ua,0)\nonumber\\
&=&d(a, 0)\cdot_{\phi}(u, 0)+d(u, 0)\cdot_{\phi}(a, 0)\\
&+&(a, 0)\cdot_{\gamma}d(u, 0)+(u, 0)\cdot_{\gamma}d(a, 0)\nonumber
\end{eqnarray}
Take $a=u$ in (1). Then
\begin{eqnarray}\label{2ud}
(z, w)&=&(z, w)\cdot_{\phi}(u, 0)+(u, 0)\cdot_{\gamma}(z, w)\nonumber\\
&=&(z+\phi(w)u+uz+\gamma(w)u, 0),
\end{eqnarray}
where $d(u, 0)=(z, w)$ for some $z\in A$ and $w\in B$. Hence $w=0$ and from (2) we obtain $az=0$ for all $a\in A$. So $z\in\mbox{ran}(A)$.

Let $d(a, 0)=(x_0, y_0)$ and $d(ua, 0)=(x_1, y_1)$ for some $x_0, x_1\in A$ and $y_0, y_1\in B$. If we replace $a$ by $ua$ in (1), then
$$
(x_1+x_1, y_1+y_1)=(x_1+\phi(y_1)u+za+\gamma(y_1)u+uaz+ux_1, 0),
$$
Hence $y_1=0$ and by (1), $y_0=0$. Therefore $d$ maps $A$ into itself.
$\hfill\square$\\

The main result of this paper is the following.

\begin{theorem} Let $d: A\times B\rightarrow A\times B$ be a mapping. Then $d\in\mbox{Der}_J(A\times_\theta^{\phi, \gamma} B)$ if and only if the following statements hold.

\emph{(i)} There exist unique Jordan derivations $d_A\in\mbox{Der}_J(A)$ and $d_B\in\mbox{Der}_J(B)$ such that
$$
d(a, b)=(d_A(a)+(2\theta-\phi-\gamma)(b)d_A(u)-\frac{1}{2}(\phi+\gamma)(d_B(b))u, d_B(b))
$$
for all $a\in A$ and $b\in B$.

\emph{(ii)} $(2\theta-\phi-\gamma)(b)(d_A(a)-d_A(u)a)=\frac{1}{2}(\phi+\gamma)( d_B(b))(a-ua)$ for all $a\in A$ and $b\in B$.

\emph{(iii)} $(\theta-\phi)(b)(\theta-\gamma)(b)d_A(u)=(\gamma-\phi)(b)(\gamma-\phi)(d_B(b))u=0$ for all $a\in A$ and $b\in B$.

\end{theorem}
{\it Proof.} For $b\in B$, let $d(0, b)=(x_1, y_1)$ and $d(u, 0)=(z, 0)$ for some $x_1\in A$, $z\in\mbox{ran}(A)$
and $y_1\in B$. Then
\begin{eqnarray*}
d(u,0)&+&d(2\theta(b)u, 0)+d(0, b^2)\\
&=&d((u, b)\cdot_{\theta}(u, b))\\
&=&d(u, b)\cdot_{\phi}(u, b)+(u, b)\cdot_\gamma d(u, b)\\
&=&d(u, 0)\cdot_{\phi}(u, 0)+d(u, 0)\cdot_{\phi}(0, b)\\
&+&d(0, b)\cdot_{\phi}(u, 0)+d(0, b)\cdot_{\phi}(0, b)\\
&+&(u, 0)\cdot_{\gamma}d(u, 0)+(u, 0)\cdot_{\gamma}d(0, b)\\
&+&(0, b)\cdot_{\gamma}d(u, 0)+(0, b)\cdot_{\gamma}d(0, b).\\
\end{eqnarray*}
So
\begin{eqnarray*}
(2\theta(b) z, 0)&=&d(2\theta(b)u, 0)\\
&=&d(u, 0)\cdot_{\phi}(0, b)+d(0, b)\cdot_{\phi}(u, 0)\\
&+&(u, 0)\cdot_{\gamma}d(0, b)+(0, b)\cdot_{\gamma}d(u, 0)\\
&=&(z, 0)\cdot_{\phi}(0, b)+(x_1, y_1)\cdot_\phi(u, 0)\\
&+&(u, 0)\cdot_{\gamma}(x_1, y_1)+(0, b)\cdot_{\gamma}(z, 0)\\
&=&(\phi(b)z+x_1+\phi(y_1)u+ux_1+\gamma(y_1)u+\gamma(b)z, 0).
\end{eqnarray*}
This shows that
\begin{eqnarray}\label{1}
x_1=(2\theta-\phi-\gamma)(b)z-ux_1-(\phi+\gamma)(y_1)u.
\end{eqnarray}
If we multiply (3) by $u$ from the left, then
$$
ux_1=-\frac{1}{2}(\phi+\gamma)(y_1)u.
$$
From this and (3), we have
$$
x_1=(2\theta-\phi-\gamma)(b)z-\frac{1}{2}(\phi+\gamma)(y_1)u.
$$
Hence
\begin{eqnarray}\label{2}
d(0, b)=((2\theta-\phi-\gamma)(b)z-\frac{1}{2}(\phi+\gamma)(y_1)u, y_1).
\end{eqnarray}
In view of Lemma 2.1, for every $a\in A$, there exists $x_0\in A$ such
that $d(a, 0)=(x_0, 0)$. This together with (4) shows that
\begin{eqnarray*}
d(a, b)&=&d(a, 0)+d(0,b)\nonumber\\
&=&(x_0+(2\theta-\phi-\gamma)(b)z\\
&-&\frac{1}{2}(\phi+\gamma)(y_1)u, y_1).\nonumber\
\end{eqnarray*}
Assume that $\pi_A: A\times B\rightarrow A$ and $\pi_B:
A\times B\rightarrow B$ be canonical projections. We define
the functions $d_A: A\rightarrow A$ and $d_B: B\rightarrow B$ by the following rules:
$$
d_A(a)=\pi_A(d(a, 0))\quad\mbox{and}\quad d_B(b)=\pi_B(d(0, b)).
$$
It is clear that these functions are Jordan derivations and
\begin{eqnarray}\label{3}
d(a, b)=(d_A(a)+(2\theta-\phi-\gamma)(b)d_A(u)-\frac{1}{2}(\phi+\gamma)(d_B(b))u, d_B(b))
\end{eqnarray}
for all $a\in A$ and $b\in B$. Hence (i) holds.

Since $d\in\mbox{Der}_J(A\times_\theta^{\phi, \gamma} B)$, for
every $a\in A$ and $b\in B$, we have
\begin{eqnarray}\label{4}
d((a, b)\cdot_\theta(a, b))=d(a, b)\cdot_\phi(a, b)+(a,
b)\cdot_\gamma d(a, b).
\end{eqnarray}
From (5) and (6), we conclude that
\begin{eqnarray}\label{5}
2\theta(b)d_A(a)&+&(2\theta-\phi-\gamma)(b^2)d_A(u)-\frac{1}{2}(\phi+\gamma)( d_B(b^2))u\nonumber\\
&=&(\phi+\gamma)(b)d_A(a)+(2\theta-\phi-\gamma)(b)d_A(u)a\nonumber\\
&+&\frac{1}{2}(\phi+\gamma)( d_B(b))(a-ua)\\
&+&(\phi+\gamma)(b)(2\theta-\phi-\gamma)(b)d_A(u)\nonumber\\
&-&\frac{1}{2}(\phi+\gamma)(b)(\phi+\gamma)( d_B(b))u\nonumber\
\end{eqnarray}
for all $a\in A$ and $b\in B$. Set $a=0$ in
(7). Then
\begin{eqnarray}\label{6}
(2\theta-\phi-\gamma)(b^2)d_A(u)&-&\frac{1}{2}(\phi+\gamma)( d_B(b^2))u\nonumber\\
&=&(\phi+\gamma)(b)(2\theta-\phi-\gamma)(b)d_A(u)\\
&-&\frac{1}{2}(\phi+\gamma)(b)(\phi+\gamma)( d_B(b))u\nonumber
\end{eqnarray}
for all $b\in B$.
Regarding (7) and (8), we infer that
\begin{eqnarray*}\label{7}
(2\theta-\phi-\gamma)(b)(d_A(a)-d_A(u)a)=\frac{1}{2}(\phi+\gamma)(d_B(b))(a-ua).
\end{eqnarray*}
That is, (ii) holds. Let us multiply (8) by $u$ from the left. Then
\begin{eqnarray*}
(\gamma-\phi)(b)(\gamma-\phi)(d_B(b))u=0
\end{eqnarray*}
for all $a\in A$ and $b\in B$. This together with (8) follows that
$$
(\theta-\phi)(b)(\theta-\gamma)(b)d_A(u)=0
$$
for all $a\in A$ and $b\in B$. Hence (iii) holds.
$\hfill\square$\\

In the sequel, $d_A$ and $d_B$ are as in Theorem 2.2.

\begin{corollary}\label{AB} Let $d$ be an element in $\mbox{Der}_J(A\times_\theta^{\phi, \gamma} B)$. Then the following statements hold.

\emph{(i)} Either $\theta=\phi=\gamma$ or $d$ maps $A$ into $\emph{\mbox{ran}}(A)$.

\emph{(ii)} If either $A$ has a unit or $A$ is semisimple, then  $\theta=\phi=\gamma$ or $d$ is zero on $A$.
\end{corollary}
{\it Proof.}In view of Theorem 2.2 (ii), for every $a, x\in A$ and $b\in B$ we have
\begin{eqnarray*}
(2\theta-\phi-\gamma)(b)x(d_A(a)-d_A(u)a)&=&\frac{1}{2}(\phi+\gamma)(d_B(b))x(a-ua)\\
&=&0.\nonumber\
\end{eqnarray*}
This implies that for every $a, x\in A$ and $b\in B$
\begin{eqnarray}\label{0}
(2\theta-\phi-\gamma)(b)xd_A(a)=0.
\end{eqnarray}
Suppose now that $d$ does not map $A$ into $\mbox{ran}(A)$. Then by (9),
$$
(2\theta-\phi-\gamma)(b)=0
$$
for all $b\in B$. Hence
$$
\theta(b)=\frac{1}{2}(\phi+\gamma)(b)
$$
for all $b\in B$. Writing $b$ by $b^2$ in the above relation, we get
$$
(\phi(b)-\gamma(b))^2=0
$$
for all $b\in B$. This implies that
$$\theta=\phi=\gamma.$$
So (i) holds. The other statement of the present result follows at once from (i).
$\hfill\square$\\

In the following, a linear mapping $d: A\times B\rightarrow A
\times B$ is called a $(\theta, \phi, \gamma)$-\emph{derivation} if
$$
d((a, b)\cdot_{\theta}(x, y))=d(a, b)\cdot_\phi(x,
y)+(a, b)\cdot_\gamma d(x, y)
$$
for all $a, x\in A$ and $b, y\in B$. The set of all these mappings is denoted by $\mbox{Der}(A\times_\theta^{\phi, \gamma} B)$.

\begin{theorem}\label{AE} Let $d$ be an element in $\mbox{Der}_J(A\times_\theta^{\phi, \gamma} B)$. Then $d\in\mbox{Der}(A\times_\theta^{\phi, \gamma} B)$ if and only if the following assertions hold.

\emph{(i)}  $d_A\in\mbox{Der}(A)$ and $d_B\in\mbox{Der}(B)$.

\emph{(ii)} $(\theta-\phi)(b)d_A(a)=(\theta-\gamma)(b)(d_A(a)-d_A(u)a)=0$ for all $a\in A$ and $b\in B$.

\emph{(iii)} $\phi d_B(b)(\phi -\gamma)(y)u=\phi d_B(b)(a-ua)=0$ for all $a\in A$ and $b, y\in B$.

\emph{(iv)} $\gamma d_B=\phi d_B$ on $B$.

Furthermore, if $A$ is a Bnach algebra without identity, then $d(a, b)=(d_A(a)+(\theta-\gamma)(b)d_A(u), d_B(b))$
for all $a\in A$ and $b\in B$.
\end{theorem}
{\it Proof.}
Let $d\in\mbox{Der}_J(A\times_\theta^{\phi, \gamma} B)$. According to Theorem 2.2, there exist $d_A\in\mbox{Der}_J(A)$ and $d_B\in\mbox{Der}_J(B)$ such that
$$
d(a, b)=(d_A(a)+(2\theta-\phi-\gamma)(b)d_A(u)-\frac{1}{2}(\phi+\gamma)(d_B(b))u, d_B(b)),
$$
for all $a\in A$ and $b\in B$.
Suppose that $d\in \mbox{Der}(A\times_\theta^{\phi, \gamma} B)$. Then for all $a,x \in A$ and $b, y\in B$
$$
d((a, b)\cdot_\theta(x, y))=d(a, b)\cdot_\phi(x, y)+(a,b)\cdot_\gamma d(x, y).
$$
So
\begin{eqnarray}\label{r}
d_A(ax)&+&\theta(b)d_A(x)+\theta(y)d_A(a)\nonumber\\
&+&(2\theta-\phi-\gamma)(by)d_A(u)-\frac{1}{2}(\phi+\gamma)(d_B(by))u\nonumber\\
&=&d_A(a)x+(2\theta-\phi-\gamma)(b)d_A(u)x-\frac{1}{2}(\phi+\gamma)(d_B(b))ux\nonumber\\
&+&\phi(y)d_A(a)+\phi(y)(2\theta-\phi-\gamma)(b)d_A(u)\\
&-&\frac{1}{2}\phi(y)(\phi+\gamma)(d_B(b))u+ad_A(x)-\frac{1}{2}(\phi+\gamma)( d_B(y))a\nonumber\\
&+&\gamma(b)d_A(x)+\gamma(b)(2\theta-\phi-\gamma)(y)d_A(u)\nonumber\\
&-&\frac{1}{2}\gamma(b)(\phi+\gamma)( d_B(y))u+\gamma d_B(y)a+\phi d_B(b)x\nonumber\
\end{eqnarray}
and
\begin{eqnarray}\label{r1}
d_B(by)=d_B(b)y+bd_B(y).
\end{eqnarray}
The relation (11) shows that $d_B$ is a derivation on $B$.
Set $b=y=0$ in (10). Then
$$
d_A(ax)=d_A(a)x+ad_A(x)
$$
for all $a,x \in A$.
Hence $d_A$ is a derivation on $A$. That is, (i) holds. Now, let $a=x=0$ in (10), we obtain
\begin{eqnarray}\label{r2}
(2\theta-\phi-\gamma)(by)d_A(u)&-&\frac{1}{2}(\phi+\gamma)(d_B(by))u\nonumber\\
&=&\phi(y)(2\theta-\phi-\gamma)(b)d_A(u)\nonumber\\
&-&\frac{1}{2}\phi(y)(\phi+\gamma)( d_B(b))u\\
&+&\gamma(b)(2\theta-\phi-\gamma)(y)d_A(u)\nonumber\\
&-&\frac{1}{2}\gamma(b)(\phi+\gamma)( d_B(y))u\nonumber\
\end{eqnarray}
Subtracting (12) from (10), we arrive at
\begin{eqnarray}\label{r3}
(\theta-\gamma)(b)(d_A(x)-d_A(u)x)&+&(\theta-\phi)(y)d_A(a)\nonumber\\
&=&(\theta-\phi)(b)d_A(u)x\nonumber\\
&-&\frac{1}{2}(\phi+\gamma)(d_B(b))ux\\
&+&\frac{1}{2}(\gamma-\phi)(d_B(y))a\nonumber\\
&+&\phi d_B(b)x\nonumber\
\end{eqnarray}
Taking $b=0$ in (13), we have
\begin{eqnarray}
(\theta-\phi)(y)d_A(a)=\frac{1}{2}(\gamma-\phi)(d_B(y))a
\end{eqnarray}
for all $a\in A$ and $y\in B$. Put $a=u$ in (14) and then multiply it by $u$ from the left. These imply that
\begin{eqnarray}
\phi d_B(b)=\gamma d_B(b)
\end{eqnarray}
for all $b\in B$. So (iv) holds. From this and (14) we infer that
\begin{eqnarray}\label{r0}
(\theta-\phi)(y)d_A(a)=0
\end{eqnarray}
for all $a\in A$ and $y\in B$. This together with (13) and (15) shows that
\begin{eqnarray}\label{r001}
(\theta-\gamma)(b)(d_A(x)-d_A(u)x)=\phi d_B(b)(x-ux)
\end{eqnarray}
for all $x\in A$ and $b\in B$. So (12) can be written as follows.
\begin{eqnarray*}\label{r5}
((\theta-\gamma)(by)&-&\phi(y)(\theta-\gamma)(b)-\gamma(b)(\theta-\gamma)(y))d_A(u)\\
&=& (\phi d_B(by)-\phi(y)\phi d_B(b)-\gamma(b)\phi d_B(y))u\nonumber\
\end{eqnarray*}
for all $b, y\in B$. Hence
$$
\phi d_R(by)-\phi(y)\phi d_B(b)-\gamma(b)\phi d_B(y)=0
$$
by Lemma 2.1. Since $d_B$ is a derivation on $B$,
\begin{eqnarray}\label{r5}
(\phi-\gamma)(b)\phi d_B(y)=0
\end{eqnarray}
for all $b, y\in B$. If $\phi d_B\neq 0$, then $\phi=\gamma$ and by (16) and (17), we get
\begin{eqnarray}\label{r6}
(\theta-\gamma)(b)(d_A(x)-d_A(u)x)=\phi d_B(b)(x-ux)=0
\end{eqnarray}
for all $x\in A$ and $b\in B$. From (16), (18) and (19) we see that the assertions (ii) and (iii) hold.

Finally, if $A$ is an algebra without identify, then by (iii) and (iv),
$$\phi d_B=\gamma d_B=0.$$
From this and (ii), we infer that
$$
d(a, b)=(d_A(a)+(\theta-\gamma)(b)d_A(u), d_B(b))
$$
for all $a\in A$ and $b\in B$.
$\hfill\square$\\

As an immediate consequence of [4], Corollary 2.3 and Theorem 2.4 we present the following result.

\begin{corollary}\label{AC} Let $A$ be a Banach algebra with identity and $B$ be a semisimple Banach algebra. If $\theta\neq\phi$, then
$\mbox{Der}_J(A\times_\theta^{\phi, \phi} B)=\mbox{Der}(A\times_\theta^{\phi, \phi} B)=\mbox{Der}(B)$.
\end{corollary}

Let us recall that a mapping $T: A\rightarrow A$ is called \emph{centralizing} if for every $a\in A$
$$
[T(a), a]\in \mbox{Z}(A),
$$
where $\mbox{Z}(A)$ is the center of $A$ and for each $a, x\in A$
$$[a, x]=ax-xa.$$
This concept can be stated for the $\theta-$Lau products $A\times_\theta B$ as follows. For nonzero multiplicative linear functionals $\eta_1, \eta_2$ on $B$, an element $d\in\mbox{Der}_J(A\times_\theta^{\phi, \gamma} B)$ is called
$(\eta_1, \eta_2)-$\emph{centralizing} if for every $a\in A$ and $b\in B$,
$$
[d(a, b), (a, b)]_{\eta_1,\eta_2}:= d(a, b)\cdot_{\eta_1}(a, b)-(a, b)\cdot_{\eta_2}d(a, b)\in\mbox{Z}(A)\times\mbox{Z}(B).
$$

\begin{theorem}\label{AD} Let $B$ be a semisimple Banach algebra. If $\theta\neq\phi$, then the only  $(\eta_1, \eta_2)-$centralizing element of $\mbox{Der}_J(A\times_{\theta}^ {\phi, \gamma}B)$ is the zero map.
\end{theorem}
{\it Proof.}
Let $d\in\mbox{Der}_J(A\times_{\theta}^ {\phi, \gamma}B)$ be $(\eta_1, \eta_2)-$centralizing. So there exist $d_A\in\mbox{Der}_J(A)$ and $d_B\in\mbox{Der}_J(B)$ such that
$$
d(a, b)=(d_A(a)+(2\theta-\phi-\gamma)(b)d_A(u)-\frac{1}{2}(\phi+\gamma)(d_B(b))u, d_B(b))
$$
for all $a\in A$ and $b\in B$. Since $d$ is  $(\eta_1, \eta_2)-$centralizing, $d_A$ and $d_B$ are centralizing on $A$ and $B$, respectively. It follows from [4, 13] that $d_B=0$ on $B$ and
$$
d_A(u)=d_A(u)u-ud_A(u)\in\mbox{Z}(A).
$$
From Lemma 2.1 we infer that
$$
d_A(u)=d_A(u)u=ud_A(u)=0.
$$
This implies that
$$
d_A(a)-ud_A(a)\in\mbox{Z}(A)
$$
and so
$$
d_A(a)-ud_A(a)=u(d_A(a)-ud_A(a))=0
$$
for all $a\in A$. Hence by Corollary 2.3 (i), we have
$$
d_A(a)=0
$$
for all $a\in A$. Therefore, $d=0$.
$\hfill\square$

\vspace{2mm}

 {\footnotesize
\noindent {\bf Mina Ghasemi}\\
Department of Mathematics,\\ Shiraz University of Technology,\\
Shiraz
71555-313, Iran\\ e-mail: mi.ghasemi@sutech.ac.ir\\
{\bf Mohammad Javad Mehdipour}\\
Department of Mathematics,\\ Shiraz University of Technology,\\
Shiraz
71555-313, Iran\\ e-mail: mehdipour@sutech.ac.ir\\\\

\footnotesize

\begin{thebibliography}{99}

\bibitem{am1} ‎ M.‎ H. Ahmadi Gandomani and M. J. Mehdipour,   \textit{Jordan, Jordan right and Jordan left derivations on convolution algebras}, Bull.‎ Iranian Math.‎ Soc.‎, {\bf 45}  (2019) 189--204.‎

\bibitem{am2} ‎ M.‎ H. Ahmadi Gandomani and M. J. Mehdipour,  \textit{Generalized derivations on some convolution algebras}, Aequationes Math.‎, {\bf 92} (2)‎  (2018), 223--241‎.‎

\bibitem{bdl} W. G. Bade, H. G. Dales and Z. A. Lykova, \textit{Algebraic and  strong  splittings  of extensions of Banach algebras}, Mem. Amer. Math. Soc., \textbf{137} (656) (1999), 39-–54.

\bibitem{b}  M. Bresar, \textit{Jordan derivations on semiprime rings}, Proc. Amer. Math. Soc., \textbf{104} (4) (1988), 1003--1006.

\bibitem{bv} M. Bresar and J. Vukman, \textit{Jordan derivations on prime rings}, Bull. Austral. Math. Soc., \textbf{37} (1988), 321--322.

\bibitem{ch}  Y. Choi, \textit{Triviality of the generalised Lau product
associated to a Banach algebra homomorphism}, Bull. Aust. Math. Soc., \textbf{94} (2016), 286–-289.

\bibitem{d} H. G. Dales, \textit{Banach Algebras and Automatic Continuity}, London Math. Soc. Monographs
24, Oxford Univ. Press, New York, 2000.

\bibitem{gm}‎ M. Ghasemi and M. J.  Mehdipour,\textit{ Derivations on Banach algebras of connected multiplicative linear functionals}, Bull.‎ Malays.‎ Math.‎ Sci.‎ Soc.‎, {\bf 44} (2021), 1727--1748‎.‎

\bibitem{JGW} J. He, J. Li, G. An and W. Huang, \textit{Characterization of $2$-local derivations and local Lie derivations of certain algebras}, Sib. Math. J., \textbf{59} (4) (2018), 721–-730.

\bibitem{h}  I. N. Herstein, \textit{Jordan derivations of prime rings}, Proc. Amer. Math. Soc., \textbf{8} (1957), 1104--1110.

\bibitem{k1} E. Kaniuth, \textit{The Bochner-Schoenberg-Eberlein property and spectral synthesis for certain Banach algebra products}, Canad. J. Math., \textbf{67} (4) (2015), 827--847.

 \bibitem{l1} A. T. Lau, \textit{Analysis on a class of Banach algebras with
applications to harmonic analysis on locally compact groups and
semigroups}, Fund. Math., \textbf{118} (3) (1983), 161--175.

\bibitem{MR}  M. Mathieu and V. Runde, \textit{Derivations mapping into the radical II}, Bull. London Math. Soc., \textbf{24} (1992), 485--487.

\bibitem{ms1} M.‎ J.‎ Mehdipour and Z.‎ Saeedi,  \textit{Derivations on group algebras of a locally compact abelian group}, Monatsh.‎ Math., {\bf 180} (2016), 595--605‎.‎

\bibitem{ms2}‎   M.‎ J.‎ Mehdipour and Z.‎ Saeedi,  \textit{Derivations on convolution algebras}, Bull.‎ Korean Math.‎ Soc.‎, {\bf 52} (2015), 123--1132‎.‎

\bibitem{sm}  M. Sangani Monfared, \textit{On certain products of Banach algebras with
applications to harmonic analysis}, Studia Math.,  \textbf{178} (3)
(2007), 277--294.

\bibitem{m} M. Sangani Monfared, \textit{Character amenability of Banach algebras}, Math. Proc. Camb. Phil. Soc., \textbf{144} (2008), 697--706.

\bibitem{sin} A. M. Sinclair, \textit{Jordan homomorphisms and derivations on semisimple Banach algebrfas}, Proc. Amer. Math. Soc., \textbf{24} (1970), 209--214.

\bibitem{w} B. Willson, \textit{Configurations and invariant nets for amenable hypergroups and related algebras}, Trans. Amer. Math. Soc., \textbf{366} (2014), 5087--5112.
\end{thebibliography}
\end{document}